\documentclass[12pt,a4paper]{article}
\usepackage{latexsym}
\usepackage{amsmath,amssymb,amsfonts,amscd,theorem}
\usepackage{enumerate}

\newcommand{\Id}{\mathrm{Id} }

\newcommand{\0}{\boldsymbol{0}}

\newcommand{\R}{\ensuremath{R} }

\newcommand{\Ad}{\operatorname{Ad}}
\newcommand{\image}{\operatorname{Im}}

\newcommand{\tr}{\operatorname{tr}}
\newcommand{\Ker}{\operatorname{Ker}}

\newcommand{\Hom}[2]{\hom(#1,#2)}

\newcommand{\TaZ}[2]{T^{\mathrm{Zar}}_{#1}(#2)}

\newcommand{\aut}{\operatorname{Aut}}
\newcommand{\Aut}[1]{\aut(#1)}

\newcommand{\SO}[1]{\mathrm{SO}_{3} (#1)}
\newcommand{\PSL}[1]{\mathrm{PSL}_2 (#1)}
\newcommand{\SL}[1]{\mathrm{SL}_2 (#1)}
\newcommand{\Sl}[1]{\mathfrak{sl}_2 (#1)}

\newcommand{\M}[1]{\mathrm{M}_3(#1)}

\newtheorem{dfn}{Definition}[section]
\newtheorem{rem}[dfn]{Remark}
\newtheorem{thm}[dfn]{Theorem}

\newtheorem{lem}[dfn]{Lemma}

\newtheorem{prop}[dfn]{Proposition}

\newtheorem{cor}[dfn]{Corollary}

\newtheorem{ex}[dfn]{Example}

\newtheorem{defn}[dfn]{Definition}

\def\Hom{\operatorname{Hom}}

\def\C{{\mathbb C}}

\def\0{\emptyset}
\def\R{{\mathbb R}}

\def\H{{\mathbb H}}

\def\N{{\mathbb N}}

\def\3{\ss}
\def\8{\infty}

\def\<{\langle}
\def\>{\rangle}

\def\co{\colon\thinspace}
   \newenvironment{proof}[1][Proof]%
   {\begin{trivlist} \item[]  {\em #1.} }%
   {\hspace*{\fill} $\Box$ \end{trivlist}}

\title{The variety of characters in $\PSL\C$}
\author{Michael Heusener and Joan Porti}
\date{}
\begin{document}

\maketitle
\begin{abstract}
We study some basic properties of the  variety of characters in
$\PSL\C$ of a finitely generated group. In particular we give an
interpretation of its points as characters of representations. We
construct 3-manifolds whose variety of characters has arbitrarily  many
components that do not lift to $\SL\C$. We also study the singular
locus of the variety of characters of a free group.\\[1ex]

\noindent {\small
}

\end{abstract}

\section{Introduction}

The varieties of representations and characters have many
applications in 3-dimensional topology and geometry. The variety of
$\SL\C$-characters has been intensively studied since the seminal
paper of Culler and Shalen \cite{CS}, but for many applications it is
more convenient to work with $\PSL\C$ instead of $\SL\C$ (see
\cite{BZ} and \cite{BMP} for instance). The purpose of this note is
to study some basic properties of the variety of characters in
$\PSL\C$. Most of the results of invariant theory that we use can be
found in any standard reference (e.g.\ \cite{KSS}, \cite{Kraft}, \cite{PV}).

Throughout this paper, $\Gamma$ will denote a finitely generated group.

\begin{defn}
   The
set of all representations of $\Gamma$ in $\PSL\C$ is denoted by \(
       R(\Gamma)
\) and it is called the \emph{variety of representations of $\Gamma$
in $\PSL\C$}.
\end{defn}

  The variety of representations $R(\Gamma)$ has a natural
structure as an affine  algebraic set over the complex numbers given as
follows: the group $\PSL\C$ is algebraic 
(see Section~\ref{sec:invariants}). Given a presentation $\Gamma=\langle
\gamma_1,\ldots,\gamma_s\mid (r_i)_{i\in I} \rangle$ we have a
natural embedding:  
\[
\begin{array}{rcl}
      R(\Gamma)&\to &\PSL\C\times\cdots\times \PSL\C\\
      \rho & \mapsto & (\rho(\gamma_1),\ldots,\rho(\gamma_s))
\end{array}
\]
and the defining equations are induced by the relations. This
structure can be easily seen to be independent of the
presentation. In fact using the isomorphism $\PSL\C\cong \SO\C$,
$R(\Gamma)$ has a structure of an affine set (see
Lemma~\ref{lem:SO3}).

The action of $\PSL\C$ on $R(\Gamma)$ by conjugation is algebraic.
The quotient $R(\Gamma)/\PSL\C$ may be not Hausdorff and it is more
convenient to consider the algebraic quotient of invariant theory,
because $\PSL\C$ is reductive.

\begin{defn}
The \emph{variety of $\PSL\C$-characters} $X(\Gamma)$ is the
quotient $R(\Gamma)/\!/\PSL\C$ of invariant theory.
\end{defn}

This definition means that $X(\Gamma)$ is an affine algebraic set
together with a regular map $t\co R(\Gamma) \to X(\Gamma)$ which
induces an isomorphism
\[
    t^{*}\co \C[X(\Gamma)]\to\C[R(\Gamma)]^{\PSL\C}
\]
   (i.e.\
the regular functions on $X(\Gamma)$ are precisely the regular
functions on $R(\Gamma)$ invariant by conjugation). We will use the
notation $R(M)=R(\pi_1M)$ and $X(M)=X(\pi_1M)$ if $M$ is a
path-connected topological space.

In this paper we study  the  basic properties of $X(\Gamma)$.

First we  explain  the name ``variety of characters": given a
representation $\rho\co\Gamma\to \PSL\C$, its character is the map
\[
\begin{array}{rcl}
      \chi_{\rho}\co\Gamma&\to&\C\\
          \gamma&\mapsto &\tr^2(\rho(\gamma))
\end{array}
\]

\begin{thm}\label{thm:PSL2characters} There is a natural bijection
between
$X(\Gamma)$ and the set of characters of representations $\rho\in
R(\Gamma)$. This bijection  maps every $t(\rho)\in X(\Gamma)$ to the
character $\chi_\rho$.
\end{thm}

In many cases the representations of $R(\Gamma)$ lift to $\SL\C$,
for instance if $\Gamma$ is a free group. In such a case,
$X(\Gamma)$ is just a quotient of the usual variety of characters
in $\SL\C$ (See Proposition~\ref{prop:naturaliso}). This quotient
is the definition already used in \cite{Bur90},
\cite{HLM1},\cite{HLM2} and \cite{Ril84} for 2-bridge knot
exteriors. The explicit computation for the figure eight knot
exterior is found in \cite{GM}.

There are cases where representations do not lift to $\SL\C$, for
instance the holonomy representation of an orientable hyperbolic
3-orbifold with 2 torsion. The next result proves that there are
manifolds with arbitrarily many components of characters that do not
lift.

\begin{thm}
\label{thm:nolifts}
   For every $n$, there exist a compact irreducible
3-manifold $M$ with $\partial M$ a 2-torus such that $X(M)$ has at
least $n$ irreducible one dimensional components  whose characters
do not lift to $\SL\C$.
\end{thm}

In Section~\ref{sec:invariants} we prove
Theorem~\ref{thm:PSL2characters}. In
Section~\ref{sec:irreducibility} we study the fiber of the
projection $t\co R(\Gamma)\to X(\Gamma)$, introducing the
different notions of irreducibility. Section~\ref{sec:lifts} is
devoted to the study of lifts
of representations and the proof of
Theorem~\ref{thm:nolifts}. In the last section we determine the
singular set of $X(\Gamma)$ when $\Gamma\cong F_n$ is the free
group of rank $n\geq 3$.

\section{Invariants of $\PSL\C$}
\label{sec:invariants}

Before proving Theorem~\ref{thm:PSL2characters} we quickly review
some basic notions of algebraic geometry and invariant theory
(that the reader may prefer to skip and go directly to the proof
in Subsection~\ref{ss:ProofTheorem}). For details see \cite{KSS},
\cite{Kraft} or \cite{PV}.

\subsection{Basic notions of invariant theory}

A closed algebraic subset $Z\subset \C^{N}$ is called  \emph{affine}.
We denote by $\C[Z]$ the ring of regular functions on $Z$. An
algebraic group $G$ that acts algebraically on $Z$ acts naturally on
$\C[Z]$ via $g f (z) := f(g^{-1} z) $. We denote by $\C[Z]^{G} $ the
ring of invariant functions, i.e.\ functions $f\in\C[Z]$ for which $g
f = f$ for all $g\in G$.

The group $G$ is called \emph{reductive} if it has the following
property: for each finite dimensional rational representation
$\rho\co G\to \mathrm{GL}(V)$ and every $G$-invariant subspace
$W\subset V$ there exist a complementary $G$-invariant subspace
$W'\subset V$, i.e.\ $V= W'\oplus W$.

If $Z$ is affine and $G$ is reductive, then the ring $\C[Z]^{G} $
is finitely generated. The affine set $Y$ such that $\C[Y]\cong
\C[Z]^{G} $ is called the \emph{algebraic quotient} and it is
denoted by  $Z /\!/\/ G$.

We shall use the following properties of reductive groups:

\begin{itemize}
\item[--] By Maschke's theorem, finite groups are reductive.

\item[--]
More generally, let $G\subset \mathrm{GL}_{n}(\C)$ be a linear
algebraic group. The group $G$ is reductive if there is a
Zariski-dense subgroup $K\subset G$ which is compact in the
classical topology. It follows that $\mathrm{GL}_{n}(\C)$,
$\mathrm{SL}_{n}(\C)$, $\textrm{O}_{n}(\C)$, $\textrm{SO}_{n}(\C)$
and $\textrm{Sp}_{n}(\C)$ are reductive.

\item[--] Let $G$ be a reductive linear algebraic group. Let $Y$ and
$Z$ be varieties  on which $G$ acts and let   $f\co X\to Y$ be a
$G$-invariant regular map. If $f^{*}\co \C[Y]\to\C[X]$ is
surjective then $f^{*}(\C[Y]^{G}) = \C[X]^{G}$ holds.
\end{itemize}

\subsection{Algebraic structure of $\PSL\C$}

The group $\PSL\C$ is algebraic, it is the quotient of $\SL\C$ by the
finite group $\{\pm \Id\}$.

It is useful to recall the isomorphism with $\SO\C$, that we
construct next. We denote by
\[
\Ad\co\PSL\C\to \Aut{\Sl\C}
\]
the adjoint action of $\PSL\C$ on its Lie algebra $\Sl\C$. The
Killing form on $\Sl\C$ is a non degenerate symmetric bilinear
form over $\C$. For each $A\in \PSL\C$,  $\Ad (A)$ preserves the
Killing form and $\det (\Ad(A)) = 1$, hence $\Ad(\PSL\C) \subseteq
\SO\C$. The following  lemma is well known from representation
theory (see for instance \cite{FultonHarris}):

\begin{lem}\label{lem:SO3}
The action of $\PSL\C$ on the Lie algebra induces an isomorphism
$\Ad\co\PSL\C\to \SO\C $.
\end{lem}

In this paper the trace will be
   abbreviated by
$\tr$,   and $\tr^{2}(A)$ stands for $(\tr(A))^{2}$.
  By direct
computation we obtain the equality
\begin{equation}\label{eqn:traceAd}
      \tr (\Ad(A))=\tr^{2}(A)-1=\tr (A^2)+1
      \quad\text{ for all } A\in \PSL\C
\end{equation}
   that will be used later.

Given $\gamma\in\Gamma$, we have a well defined function
\[
\begin{array}{rcl}
      \tau_{\gamma}\co R(\Gamma)&\to&\C\\
                                       \rho&\mapsto &\tr^{2}(\rho(\gamma))
\end{array}
\]
Since it is invariant by conjugation, it induces a function
\[
J_{\gamma}\co X(\Gamma) \to\C.
\]

\subsection{Proof of Theorem~\ref{thm:PSL2characters}}
\label{ss:ProofTheorem}

  Theorem~\ref{thm:PSL2characters} is a consequence of:

\begin{prop}\label{prop:generators} The ring of invariant functions
$\C[R(\Gamma)]^{\PSL\C}$ is generated by
the functions $\tau_{\gamma}$, with $\gamma\in\Gamma$.
\end{prop}

\begin{proof}
There is a surjection $\psi\co F_{n}\to \Gamma$ where  $F_n$ is a
free group of rank $n\in\N$. We obtain an inclusion 
$\psi^{*}\co R(\Gamma) \subset R(F_{n})$. This inclusion induces a surjection
$\psi_{*} \co\C[R(F_{n})]\to \C[R(\Gamma)]$. Now, $\PSL\C$ is
reductive and acts regularly by conjugation on the representation
varieties. Hence we obtain a surjection
\[
\psi_{*}\co \C[R(F_{n})]^{\PSL\C}\to \C[R(\Gamma)]^{\PSL\C}
\]
and it is sufficient to prove the proposition for $\Gamma =
F_{n}$ since $\psi_{*} (\tau_{\gamma}) = \tau_{\psi(\gamma)}$.

Using Lemma~\ref{lem:SO3} and (\ref{eqn:traceAd}), we have to prove
that $\C[R(F_n)]^{\SO\C}$ is generated by the trace
functions on elements of $F_n$. Equivalently, we claim that
\[
\C[\SO\C\times\cdots\times\SO\C]^{\SO\C}
\]
is generated by traces of products of matrices  and their transposes.

Let $\M\C$ denote the algebra of $3\times 3$ matrices with complex
coefficients. The group $\PSL\C\cong \SO\C$ acts on the product
$\M\C\times \cdots \times\M\C$ diagonally by conjugation.
   A theorem of Aslaksen, Tan and Zhu (see \cite{ATZ})
     states that
the algebra of invariant functions
\[
\C[\M\C\times \cdots \times \M\C]^{\SO\C}
\]
is generated by the traces of products of matrices and their
transposes. 
Thus the proof of the proposition reduces  to show that
we have a natural surjection
\[
\C[\M\C\times \cdots \times \M\C]^{\SO\C} \to
\C[\SO\C\times \cdots \times \SO\C]^{\SO\C}\,.
\]
Since $\SO\C\times \cdots \times \SO\C\subset \M\C\times \cdots \times
\M\C$
is a closed subvariety we obtain a natural surjection
\[
\C[\M\C\times \cdots \times \M\C]\to \C[\SO\C\times \cdots\times \SO\C]
\]
which is of course $\SO\C$-invariant.
Using the fact that  $\SO\C$ is reductive gives the surjection
$\C[\M\C\times \cdots \times \M\C]^{\SO\C} \to
\C[\SO\C\times \cdots \times \SO\C]^{\SO\C}$.
\end{proof}

Since $\C[X(\Gamma)]=\C[R(\Gamma)]^{\SO\C}$ is finitely generated, we
also obtain:

\begin{cor}\label{cor:embedding} There are finitely many elements
$\gamma_{1},\ldots,\gamma_{N}$ in $\Gamma$ such that
$J_{\gamma_1}\times\cdots\times J_{\gamma_N}\co X(M)\to\C^N$ is an
embedding and  its image is
a closed algebraic set.
\end{cor}

\subsection{Other invariant functions}
\label{ssec:others} 

 There are other natural
functions to consider. Let $\Gamma^{2}$ be 
the subgroup of $\Gamma$
generated by the squares $\gamma^{2}$ of all elements $\gamma$ of
$\Gamma$. It is well known that we have an exact sequence:
\[
1\to\Gamma^2\to\Gamma\to H_1(\Gamma,{C_2})\to 1,
\]
 where $C_2=\{\pm 1\}$ is the group with $2$ elements.
   For instance, if $\Gamma$ is a finite group of odd order,
  then $\Gamma^{2}=\Gamma$.
In general, if $ \gamma,\mu \in\Gamma$ the commutator
$[\gamma,\mu]=\gamma\mu\gamma^{-1}\mu^{-1}=
(\gamma\mu)^{2}(\mu^{-1}\gamma^{-1}\mu)^{2}\mu^{-2}$ is in $\Gamma^{2}$
and hence $\Gamma^{2}$ contains the commutator group
$\Gamma' = [\Gamma,\Gamma]$. Notice
that
$$\Gamma^2=\bigcap_{\epsilon\in H^1(\Gamma,C_2)} \Ker(\epsilon)$$
where
$H^1(\Gamma,C_2)=\Hom(\Gamma,C_{2})$. Let $R(\Gamma,\SL\C)$ denote
the variety of representations of $\Gamma$ in $\SL\C$. The
cohomology group $H^1(\Gamma,{C_2})$ acts on this variety of
representations as follows: an homomorphism
$\epsilon\co\Gamma\to{C_2}=\{\pm 1\}$ maps the representation
$\rho\in R(\Gamma,\SL\C)$ to the product of representations
$\epsilon\cdot \rho$ (which maps $\gamma\in\Gamma$ to $
{\epsilon(\gamma)}\cdot \rho(\gamma)$).

\subsubsection{Invariant functions for the free group}

Let $F$ be a finitely generated free group.
For $\gamma\in F^{2}$ and $\rho\in R(F)$,
$\tr(\rho(\gamma))$ is well defined since
the representation $\rho\co F \to \PSL\C$ lifts to
$\tilde{\rho}\co F \to \SL\C$ and for $\gamma\in F^{2}$ the trace
$\tr(\tilde{\rho}(\gamma))$ depends only on $\gamma$.
Note that two lifts $\tilde{\rho}_{1}$ and $\tilde{\rho}_{2}$ of
$\rho$ differ by a homomorphism $\epsilon\in H^1(F,C_2) $ and that
$F^{2}\subset \Ker(\epsilon)$ for each
$\epsilon\in H^{1}(F,C_{2})$.

\begin{prop}
\label{prop:sigmafunction}
Let $F$ be a free group.
  For every $k$-tuple $\gamma_1,\ldots,\gamma_k\in F$ such that the product
  $\gamma_1\cdots\gamma_k\in F^2$, the function
\[
\begin{array}{rcl}
     \sigma_{\gamma_1,\ldots,\gamma_k}\co R(F)&\to&\C \\
     \rho&\mapsto&\tr(\tilde\rho(\gamma_1))\cdots
     \tr(\tilde\rho(\gamma_k))
\end{array}
\]
is regular (i.e.\
$\sigma_{\gamma_1,\ldots,\gamma_k}\in\C[R(F)]$).
Here,  $\tilde\rho\co F \to\SL\C$ denotes a lift of $\rho$.
\end{prop}

In order to prove this proposition we shall use the following:
\begin{lem}
\label{lem:RFn}
  Let $F_n$ be the free group of rank $n$. We have a
natural isomorphism
\[
  R(F_n,\SL\C)/\! /H^1(F_n,{C_2})\cong R(F_n)\,.
\]
\end{lem}

\begin{proof}
     Since $R(F_n,\SL\C)\cong \SL\C^n$, $R(F_n)\cong \PSL\C^n$
     and $\SL\C/{C_2}\cong \PSL\C $, we have the lemma.
\end{proof}

\begin{proof}[Proof of Proposition~\ref{prop:sigmafunction}]
For a free group $F$ and $\gamma_1,\ldots,\gamma_k\in F$, the
function $\tilde \sigma \co R(F,\SL\C)\to \C$ given by
$\tilde\sigma(\rho)=\tr(\rho(\gamma_1))\cdots \tr(\rho(\gamma_k))$
is regular. Moreover, we have $\tilde\sigma(\epsilon\cdot\rho) =
\epsilon(\gamma_1\cdots\gamma_k) \tilde\sigma(\rho)$. Since the
product $\gamma_1\cdots\gamma_k\in F^{2}$ we get that
$\tilde\sigma\in \C[R(F_n,\SL\C)]^{H^1(F_n,{C_2})}$ is an
invariant regular function on the $\SL\C$ representation variety.
By Lemma~\ref{lem:RFn}, this function factors through $R(F)$ and
gives the regular function
$\sigma_{\gamma_1,\ldots,\gamma_k}\in\C[R(F)] $.
\end{proof}

\begin{ex}
\label{ex:sigmagammaeta} Given $\gamma,\eta\in F$, by
Proposition~\ref{prop:sigmafunction},
$\sigma_{\gamma,\eta,\gamma\eta}\in \C[R(F)]$, thus by
Proposition~\ref{prop:generators}, $\sigma_{\gamma,\eta,\gamma\eta}$
is a polynomial on the functions $\tau$.
\end{ex}

To compute explicitly the polynomial of
Example~\ref{ex:sigmagammaeta}, we recall some identities of traces
in $\SL\C$:
\[
\tr(AB)=\tr (BA)\quad\textrm{and}\quad \tr (A)=\tr (A^{-1})\qquad
\forall A,B\in \SL\C\,.
\]
In addition, we have the fundamental identity:
\begin{equation}\label{eqn:sl2fundamental}
      \tr (AB)+\tr (A^{-1}B)=\tr (A)\,\tr (B)
      \qquad \forall A,B\in
\SL\C .
\end{equation}
This identity can be deduced from $A^2-(\tr A)A+\Id=0$
   multiplying by $A^{-1}B$  and taking traces.
Taking the square of $\tr(AB^{-1})=\tr (A)\,\tr
(B)-\tr(AB)$ we deduce:
\begin{equation*}
2 \tr(A)\tr(B)\tr(AB) = \tr^{2}(A) \, \tr^{2}(B)+\tr^{2}(AB) -
\tr^{2}(AB^{-1})\,.
\end{equation*}
Thus
\begin{equation}
     \label{eqn:sigmaABAB}
\sigma_{\gamma,\eta,\gamma\eta}=
\frac12(\tau_{\gamma}\tau_{\eta}+\tau_{\gamma\eta}-\tau_{\gamma\eta^{-1}}).
\end{equation}

\begin{ex} For every $\gamma,\mu\in F$, the commutator
$[\gamma,\mu]=\gamma\mu\gamma^{-1}\mu^{-1}$ belongs to $F^2$ and
therefore $\sigma_{[\gamma,\mu]}\in\C[R(F)]$. Using the the same
method as for Equation~(\ref{eqn:sigmaABAB}) one can find:
\begin{equation}\label{ex:traces}
\sigma_{[\gamma,\eta]}= \tau_{\gamma} + \tau_{\eta} + \frac 1 2
\tau_{\gamma\eta} + \frac 1 2 \tau_{\gamma\eta^{-1}} - \frac 1 2
\tau_{\gamma}\tau_{\eta} -2 \,.
\end{equation}
\end{ex}

\subsubsection{Invariant functions for other groups}

Let $\Gamma$ be a finitely generated group, $F$   a free group and
 $\psi\co F\to\Gamma$ a surjection. It induces another
surjection $\psi_{*}\co \C[ R(F)]\to \C[R(\Gamma)]$, $\psi_{*}f
(\rho) = f (\rho\circ \psi)$. Hence we obtain for all
$\eta_{1},\ldots,\eta_{k}\in F$ such that the product
$\eta_{1}\cdots\eta_{k}\in F^{2}$ a regular function
$\psi_{*}\sigma_{\eta_{1},\ldots,\eta_{k}}\in \C[R(\Gamma)]$.
  Note that the functions
$\psi_{*}\sigma_{\eta_1}$ and $\psi_{*}\sigma_{\eta_2}$ might be
different even if $\psi(\eta_1)=\psi(\eta_2)$ in $\Gamma$. This
reflects the fact that in general not every representation
$\rho\co\Gamma\to\PSL\C$ lifts to $\SL\C$.
\begin{ex}
Let $\psi\co F\to\Gamma$ be the canonical projection
where $F=\langle x,y \mid - \rangle$ and
$\Gamma = \langle x,y \mid [x,y] =1 \rangle$.
We consider the representation
$\rho\co\Gamma\to\PSL\C$ given by
$ \rho (x) = \pm A_{x}$ and $\rho (y) = \pm A_{y}$ where
$$A_{x} =\begin{pmatrix} i & 0 \\ 0& -i \end{pmatrix} \text{ and }
        A_{y} = \begin{pmatrix} 0 & 1 \\ -1& 0 \end{pmatrix}\,.$$
We obtain $\tr ([A_{x},A_{y}]) = -2$ and hence
$\psi_{*}\sigma_{[x,y]}(\rho) = -2$. On the other hand
we have $[x,y] =1$ in $\Gamma$ and
 $\psi_{*}\sigma_{1}= 2$ is a constant function.
\end{ex}

If the representation $\rho\in R(\Gamma)$ admits a lift
$\tilde\rho\co\Gamma\to\SL\C$ then
\begin{equation}\label{equ:sigma}
\psi_{*}\sigma_{\eta_{1},\ldots,\eta_{k}}(\rho)=
\tr(\tilde\rho(\psi(\eta_{1})))\cdots
\tr(\tilde\rho(\psi(\eta_{k})))
\end{equation}
only depends on the elements $\psi(\eta_1),\ldots,
\psi(\eta_k)\in\Gamma$.

\section{Irreducibility}
\label{sec:irreducibility}

To study the fiber of the map $t\co R(\Gamma)\to X(\Gamma)$ we shall
consider two different notions of irreducibility for $\rho\in
R(\Gamma)$, the usual one as a
representation in $\PSL\C$ and the so
called $\Ad$-irreducibility, for the three dimensional representation
$\Ad\circ\rho\co \Gamma\to SO_3(\C)$.

\subsection{Irreducible representations}

\begin{defn} A representation $\rho\in R(\Gamma)$ is called
\emph{reducible} if $\rho(\Gamma)$ preserves a
          point of $P^1(\C)$. Otherwise it is called irreducible.
          A character $\chi\co\Gamma\to\C$ is called \emph{reducible} if
it is the character of a reducible representation.
\end{defn}

\begin{rem}\label{rem:reducible}
Up to conjugation, the image of a reducible representation is
contained in the set of upper-triangular matrices
$\left(\begin{smallmatrix}
* & *\\ 0 & *\end{smallmatrix} \right)$.
\end{rem}

We shall require the following well known lemma (see
\cite[\S~4.3]{Beardon}).

\begin{lem}
\label{lem:fixedpoint}
  Two non-trivial elements $g,h\in\PSL\C$ have a common
fixed point in $P^{1}(\C)$ if and only if $\tr([g,h])=2$. In
addition, this fixed point is unique if $[g,h]$ is not the
identity.
\end{lem}

Irreducibility is a property that can be detected from characters:

\begin{lem}
\label{lem:tracered}
A representation $\rho\in R(\Gamma)$ is reducible
iff $\tr([\rho(\gamma),\rho(\eta)])=2$ for all elements
$\gamma,\eta$ in $\Gamma$.
\end{lem}
\begin{proof}
If $\rho$ is reducible then all the $\rho(\gamma)$ have a common
fixed point and Lemma~\ref{lem:fixedpoint} gives the result.

Assume now that $\tr([\rho(\gamma),\rho(\eta)])=2$ for all elements
$\gamma,\eta$ in $\Gamma$.

\smallskip

\emph{Case 1:} There are two elements $\gamma$ and $\eta$ in
$\Gamma$ such that  $[\rho(\gamma),\rho(\eta)]$ is not the
identity. Then $A=[\rho(\gamma),\rho(\eta)]$  is a non-trivial
parabolic element in the image of $\Gamma$.  For any
$\mu\in\Gamma$, either $\rho(\mu)$ commutes with $A$ or $[\rho(\mu),A]$ is non-trivial. 
The former possibility implies
that $\rho(\mu)$ fixes the unique fixed point of $A$, the latter
too by Lemma~\ref{lem:fixedpoint}.

\smallskip

\emph{Case 2:} The image of $\rho$ is an abelian group. Abelian
subgroups of $\PSL\C$ are well-known: either they have a global fixed
point in $P^{1}(\C)$ or they are conjugated to the group with four
elements generated by $\pm\left(\begin{smallmatrix} 0 & 1\\ -1 &
0\end{smallmatrix} \right)$ and $\pm \left(\begin{smallmatrix} i &
0\\ 0 & -i\end{smallmatrix} \right)$. 
Since the commutator of these
two generators is $\left(\begin{smallmatrix} -1 & 0\\ 0 &
-1\end{smallmatrix} \right)$, this possibility does not occur.
\end{proof}

\begin{dfn}
  \label{def:Klein}
  A non-cyclic  abelian subgroup of $\PSL\C$ with four elements is
  called \emph{ Klein's 4-group}. Such a group is realized by
  rotations about three orthogonal geodesics  and it is conjugated
  to the one generated by 
  $\pm\left(\begin{smallmatrix} 0 & 1\\ -1 & 0\end{smallmatrix} \right)$ and 
  $\pm \left(\begin{smallmatrix} i & 0\\ 0 & -i\end{smallmatrix} \right)$.
\end{dfn}

Let $R^{red}(\Gamma)$ denote the set of reducible representations
and $X^{red}(\Gamma)=t(R^{red}(\Gamma))$.
Let $F$ be a free group and let $\psi\co F\to\Gamma$ be surjective. Lemma~\ref{lem:tracered}
implies that
$$R^{red}(\Gamma) = \{\rho\in R(\Gamma)\mid
\psi_{*}\sigma_{[\gamma,\eta]}(\rho)=2 \quad \forall \gamma,\eta\in F \}$$
is a Zariski closed subset
invariant by conjugation. Thus, by invariant theory we have:

\begin{cor}
\label{cor:red}
   The set $X^{red}(\Gamma)$ is Zariski closed and
$R^{red}(\Gamma)=t^{-1}(X^{red}(\Gamma))$.
\end{cor}

\begin{rem}
\label{rem:rhodiagonal} Every reducible character $\chi$ is the
character of a diagonal representation, because if
$\rho(\gamma)=
        \pm\left(\!\begin{smallmatrix}a_{\gamma}&b_{\gamma}\\
                                      0&
c_{\gamma}\end{smallmatrix}\!\right)$
is a representation, then
$\rho'(\gamma)=\pm\left(\!\begin{smallmatrix}a_{\gamma}&0\\ 0&
c_{\gamma}\end{smallmatrix}\!\right)$
is also a  representation with
$\chi_{\rho}=\chi_{\rho'}$.
\end{rem}

\subsection{$\Ad$-irreducibility}

\begin{defn} A representation $\rho\in R(\Gamma)$ is
\emph{$\Ad$-reducible} if $\Sl\C$ has a proper
          invariant subspace by the action of
          $\Ad\circ\rho$. Otherwise it is $\Ad$-irreducible.
\end{defn}

Let $\H^3$  denote the three-dimensional hyperbolic space and
$\partial_{\infty}\H^3$ its ideal boundary. We use the isomorphism
$\operatorname{Isom}^+(\H^3)\cong \PSL\C$ and the natural
identification $\partial_{\infty}\H^3\cong  P^1(\C)$.

\begin{lem} A representation $\rho\co\Gamma\to\PSL\C$ is
$\Ad$-reducible if and only if 
$\rho(\Gamma)$ preserves either
a point  in $\partial_{\infty}\H^3$ or a geodesic
in $\H^{3}$.
\end{lem}

\begin{proof} Let $V$ be a proper subspace of $\Sl\C$ invariant by
$\Ad\circ\rho(\Gamma)$. Up to taking $V^{\bot}$, we may assume
$\dim V=1$, because the Killing form   is not degenerate. We have
then two possibilities: either the Killing form restricted to $V$
vanishes or not. In the first case $V$ consists of parabolic
Killing fields, in particular the 1-parameter group $\{\exp(v)\mid
v\in V\}\cong \C$ is parabolic and fixes a unique point at
infinity, that has to be fixed also by $\rho$. In the second case,
when the Killing form restricted to $V$ does not vanish,   the
1-parameter group $\{\exp(v)\mid v\in V\}\cong \C^{*}$ is a subgroup
of index two in the group of isometries which preserve a geodesic in $\H^3$.
This geodesic has to be preserved by the representation.
 Conversely, if
a representation preserves a point in $\partial_{\infty}\H^3$ or a
geodesic, the previous argument shows how to construct an
invariant subspace of $\Sl\C$.
\end{proof}

\begin{cor} Reducible representations are also $\Ad$-reducible.
\end{cor}

\begin{rem}
A representation $\Ad$-reducible but not reducible is a
${C_2}$-extension of an abelian one that fixes an oriented geodesic.
Thus it preserves an unoriented geodesic
\end{rem}

We call a representation $\rho\in R(\Gamma)$ \emph{abelian}
respectively \emph{metabelian} if its image is an abelian
respectively metabelian subgroup of $\PSL\C$

\begin{lem} \label{lem:Ad-red}
   A representation $\rho\in R(\Gamma)$ is $\Ad$-reducible
iff it is metabelian.
\end{lem}
\begin{proof}
   If  $\rho$ is $\Ad$-reducible then its image is contained in the
   stabilizer of either a point in $P^1(\C)$ or a geodesic in $\mathbb
   H^3$. Those stabilizers are metabelian, since they are respectively
   the group of affine transformations of $\C$ and the semidirect
   product $\C^{*}\rtimes C_{2}$.

   Now assume that $\rho(\Gamma)\subset\PSL\C$ is a metabelian
   subgroup. We use the fact that an abelian subgroup of $\PSL\C$
    preserves a unique point of $P^1(\C) $, a unique geodesic or it is
Klein's 4-group
   (Def.~\ref{def:Klein}). If $\rho([\Gamma,\Gamma])$ is trivial then
    $\rho$ is $\Ad$-reducible by this fact. If
   $\rho([\Gamma,\Gamma])$is not trivial, then we look at those
   \emph{unique} invariant objects: the \emph{unique} point in
$P^1(\C) $, the \emph{unique}
   geodesic, or the \emph{unique} three geodesics if it is Klein's 4-group.
   Since $[\Gamma,\Gamma]$ is normal in $\Gamma$, uniqueness implies that
   $\rho(\Gamma)$ preserves the same objects, hence $\rho$ is
   $\Ad$-reducible.
\end{proof}

\begin{lem}
   The set of characters of $\Ad$-reducible representations is Zariski
   closed.
\end{lem}

\begin{proof}
Lemma~\ref{lem:Ad-red} gives that the set of $\Ad$-reducible
representations is $$ R^{\Ad-red}=\{ \rho\in R(\Gamma)\mid \rho(c) =
\pm\Id
         \quad \forall c \in \Gamma''\}$$
where $\Gamma''$ denotes the second commutator group of $\Gamma$.
This is a closed subset of $R(\Gamma)$ invariant under
conjugation. Hence we have $X^{\Ad-red}(\Gamma)=t(R^{\Ad-red})$ is
a closed subset of $X(\Gamma)$.
\end{proof}

\begin{rem} The image of an $\Ad$-reducible representation is
elementary, but elementary groups also include groups that fix a
point in $\mathbb H^3$.
\end{rem}

\subsection{The fibers of $t\co R(\Gamma)\to X(\Gamma)$}

\begin{lem}\label{lem:fiberirred}
The fiber of an irreducible character consists of a single closed
orbit (i.e.\ two irreducible representations have the same character
iff they are conjugate).
\end{lem}

\begin{proof} Let $\rho_1,\rho_2\in R(\Gamma)$ be two
irreducible representations with $\chi_{\rho_1}=\chi_{\rho_2}$.

We assume first that each $\rho_i$ is irreducible but $\Ad$-reducible.
Thus each $\rho_i$ preserves a geodesic $l$, that we may assume to be
the same after conjugation. The action of $\rho_i(\gamma)$ on $l$ is
determined by the value of $\chi_{\rho_i}(\gamma)$, except in the
case $\chi_{\rho_i}(\gamma)=0$, which means that $\rho_i(\gamma)$ is
a rotation through angle $\pi$, but
    it can be either about $\gamma$ or about an axis perpendicular to
$\gamma$.
Thus if there exists an element $\gamma_0\in \Gamma$ with
$\chi_{\rho_i}(\gamma_0)\neq 4,0$ (i.e.\ $\rho_i(\gamma_0)$ is
either a loxodromic element or a rotation of angle $\neq\pi$) then
$\forall\gamma\in\Gamma$ the action of $\rho_i(\gamma)$ on the
geodesic $l$ is determined by $ \chi_{\rho_i}(\gamma)$ and
$\chi_{\rho_i}(\gamma\gamma_0)$. In particular $\rho_i$ is unique
up to conjugation. The exceptional case occurs when
$\chi_{\rho_i}(\gamma)=0$ or $4$ for every $\gamma\in\Gamma$. In
this special case, $\rho_i$ is necessarily a representation into
Klein's 4-group. The lemma is also clear in this case.

When $\rho_i$ are $\Ad$-irreducible, we can   assume that $\Gamma$
is a free group. Thus we can lift $\rho_i$ to  \(
\tilde\rho_i\co\Gamma\to\SL\C \).
   By Example~\ref{ex:sigmagammaeta}, for every
pair $\gamma,\gamma'\in \Gamma$ we obtain a regular function
$\sigma_{\gamma,\gamma',\gamma\gamma'}\co X(\Gamma)\to\C$, given
by
$$\sigma_{\gamma,\gamma',\gamma\gamma'}(\chi_{\rho}) =
   \tr\tilde\rho(\gamma\gamma')\,\tr\tilde\rho(\gamma)\,
    \tr\tilde\rho(\gamma')$$
    where $\tilde\rho\co\Gamma\to\SL\C$ is any lift of $\rho$.
Thus:
\begin{equation}
\label{eqn:traceslifts}
\tr\tilde\rho_1(\gamma\gamma')\,
\tr\tilde\rho_1(\gamma)\,
\tr\tilde\rho_1(\gamma') =
\tr\tilde\rho_2(\gamma\gamma')\,
\tr\tilde\rho_2(\gamma)\,
\tr\tilde\rho_2(\gamma')\, .
\end{equation}

We define $\epsilon\co\Gamma\to{C_2}=\{\pm 1\}$ by the formula:
\[
\phantom{aa}\qquad \tr\tilde\rho_1(\gamma)= {\epsilon(\gamma)}
\tr\tilde\rho_2(\gamma), \qquad \forall\gamma\in\Gamma\textrm{
such that } \chi_{\rho_1}(\gamma)\neq 0.
\]
When $\chi_{\rho_1}(\gamma)=0$, since we assume that $\rho_i$ is
$\Ad$-irreducible, we can find $\gamma_0\in \Gamma$ with
$\chi_{\rho_i}(\gamma_0)\neq 0$ and
$\chi_{\rho_i}(\gamma\gamma_0)\neq 0$. In this case we define
$\epsilon (\gamma)=\epsilon(\gamma_0)\cdot \epsilon(\gamma
\gamma_0)$.

By (\ref{eqn:traceslifts}), $\epsilon$ is a morphism.
   Hence $\tilde\rho_1$
and $ {\epsilon}\cdot \tilde\rho_2$ are irreducible
representations in $\SL\C$ with the same character. By \cite{CS}
they are conjugate.
\end{proof}

\begin{prop}
\label{lem:stabilizers}
\begin{itemize}
\item[(i)] A character $\chi$ is irreducible iff $\PSL\C$ acts
transitively
on the fiber and with finite stabilizer.
\item[(ii)] A character is $\Ad$-irreducible iff $\PSL\C$ acts faithfully
on the fiber.
\end{itemize}
\end{prop}

\begin{proof} (i) By
Lemma~\ref{lem:fiberirred}, if $\chi$ is irreducible then $\PSL\C$
acts transitively on $t^{-1}(\chi)$. Assume now that the stabilizer
is infinite: i.e.\ there exists nontrivial $A\in \PSL\C$ of order
$\geq 3$ (possibly infinite) and $\rho$ in the fiber such that $A$
commutes with $\rho$. If $A$ is parabolic, then it has a fixed point
in $P^1(\C)$ and therefore $\rho$ fixes this point. Otherwise $A$ has
an invariant geodesic; since $A$ has order $\geq 3$, $\rho$ preserves
the oriented geodesic, and therefore $\rho$ is also reducible.

Assume the character is reducible, then it has a diagonal
representation   $\rho$ on the fiber (Rem.~\ref{rem:rhodiagonal}),
and therefore the group of diagonal matrices stabilizes it. Thus the
stabilizer is infinite.

\smallskip

(ii) Assume $\PSL\C$ does not act 
faithfully on the fiber, i.e.\ there
exists nontrivial $A\in \PSL\C$ and $\rho$ in the fiber such that $A$
commutes with $\rho$. If $A$ is parabolic, then $\rho$ fixes a point
in $P^1(\C)$ by the previous argument. Otherwise $A$ has an
invariant geodesic, and by commutativity, $\rho$ must preserve this
geodesic. In both cases, $\rho$ is $\Ad$-reducible.

If the character is irreducible but $\Ad$-reducible, then it
preserves a geodesic, and the rotation through angle $\pi$ about this
geodesic commutes with $\rho$. 
Hence the stabilizer is nontrivial.
\end{proof}

\begin{rem} The projection $t\co R(\Gamma)\to X(\Gamma)$ induces a
bijection between irreducible components.
\end{rem}

A priori 
$R(\Gamma)$ could have more components than $X(\Gamma)$, but
the number of components is the same, because $\PSL\C$ is
irreducible.

  From Corollary~\ref{cor:red} and Proposition~\ref{lem:stabilizers} we
deduce:

\begin{cor}
Let $\rho\in R(\Gamma)$ be an irreducible representation. Let $R_0$
denote an irreducible component of $R(\Gamma)$ that contains $\rho$
and let $X_0$ denote the corresponding irreducible component of
$X(\Gamma)$. Then
\[
\dim R_0=\dim X_0+3.
\]
\end{cor}

\section{Lifts of representations to $\SL\C$}
\label{sec:lifts}

Let $\overline R(\Gamma)\subset R(\Gamma)$ denote the set of
representations $\rho\in R(\Gamma)$ that lift to $\SL\C$. According
to \cite[Thm.~4.1]{Culler} $\overline R(\Gamma)$ is a union of
connected components of $R(\Gamma)$. In particular  $\overline R
(\Gamma) $ is a Zariski-closed algebraic subset of $R(\Gamma)$, since
irreducible complex varieties are connected in the $\C$-topology 
\cite[VII, \S 2]{Sha}.
Moreover, $\overline R (\Gamma)$ is invariant under conjugation and
hence the algebraic quotient $$ \overline X(\Gamma)=\overline
R(\Gamma)/\! /\PSL\C$$ is a well defined closed subset of
$X(\Gamma)$.

In many cases, $\overline X(\Gamma)=X(\Gamma)$. For instance this is
clear when $\Gamma$ is a free group. It is also true if
$H^2(\Gamma,{C_2})=0$ by the following remark (see \cite{GM} or
\cite{Culler}).

\begin{rem}
Let $\rho\co\Gamma\to\PSL\C$ be a representation. There is a
second Stiefel-Whitney class $w_{2}(\rho)\in H^{2}(\Gamma,{C_2})$
which is exactly the obstruction for the existence of a  lift
$\overline \rho\co\Gamma\to\SL\C$.
\end{rem}

\subsection{Properties of $\overline X(\Gamma)$}

Let $R(\Gamma,\SL\C)$ and $X(\Gamma,\SL\C)$ denote the variety of
representations and characters in $\SL\C$. The ring
$\C[R(\Gamma,\SL\C)]^{\SL\C}$ is generated by the trace functions
$\tilde\tau_{\gamma}\co R(\Gamma,\SL\C)\to\C$,
$\tilde\tau_{\gamma}(\rho)=\tr(\rho(\gamma))$. The function
induced by $\tilde\tau_{\gamma}$ is denoted by $I_{\gamma}\co
X(\Gamma)\to\C$, therefore $\C[X(\Gamma)]$ is finitely generated
by the functions $I_{\gamma}$, $\gamma\in\Gamma$ \cite{CS}.

Elements of the cohomology group $  H^1( \Gamma,{C_2})$ are 
homomorphisms
$\theta\co\Gamma\to{C_2}=\{\pm 1\}$ that act on representations by
multiplication. The action of $\epsilon \in H^1( \Gamma,C_2)$ on
$I_{\gamma}$ is given by: $ \epsilon \cdot I_{\gamma} = \epsilon
(\gamma) I_{\gamma}$. Since $ H^1( \Gamma,{C_2})$ is finite, it
is reductive and we may take the quotient of invariant theory.

Let $F$ be a finitely generated free group and $\psi\co F\to\Gamma$
be a surjection. We fix a $k$-tuple $\gamma_{1},\ldots,\gamma_{k}\in
\Gamma$
such that the product $\gamma_{1}\cdots\gamma_{k}\in\Gamma^{2}$.
 Moreover, we choose $\eta_{i}\in F$ such that
$\psi(\eta_{i})=\gamma_{i}$ and such that the product
$\eta_{1}\cdots\eta_{k}\in F^{2}$. The function
$\psi^{*}\sigma_{\eta_1,\ldots,\eta_k}\in \C[\overline R(\Gamma)]$
is invariant under conjugation and gives us a function
$\psi^{*}\sigma_{\eta_1,\ldots,\eta_k}\in\C[\overline X(\Gamma)]$.
By Equation~(\ref{equ:sigma}) we have
$\psi^{*}\sigma_{\eta_1,\ldots,\eta_k} (\chi) =
\tilde\chi(\gamma_{1})\cdots\tilde\chi(\gamma_{k})$ where
$\tilde\chi\in X(\Gamma,\SL\C)$ is a character such that
$\pi(\tilde\chi)= \chi$. Note that $\pi\co X(\Gamma,\SL\C)\to
\overline X(\Gamma)$ is surjective. The function
\begin{equation}
\label{eqn:Sigma}
 \Sigma_{\gamma_1,\ldots,\gamma_k} :=
\phi^{*}\sigma_{\eta_1,\ldots,\eta_k} \in\C[\overline X(\Gamma)]
\end{equation}
depends only  on the elements $\gamma_{i}\in\Gamma$.

\begin{prop}
\label{prop:naturaliso}
 There is a natural isomorphism:
\[
X( \Gamma,\SL\C)/\! /H^1( \Gamma,C_2)\cong \overline X(\Gamma).
\]
\end{prop}

\begin{proof}
     Composition with the projection $\SL\C\to\PSL\C$ induces a
     surjection
     \[
         \pi\co X(\Gamma,\SL\C) \to \overline X(\Gamma),
     \]
     which is easily seen to be algebraic and is given by
     $\pi(\chi) = \chi^{2}$.
     At the level of function
     rings it induces an injection
     \[
     \pi^*\co  \C[\overline X(\Gamma)]\hookrightarrow \C[X(\Gamma,\SL\C)].
     \]
     We have $\pi^{*}f (\chi)= f(\chi^{2})$ for
     $f\in\C[\overline X(\Gamma)]$ and $\chi\in X(\Gamma,\SL\C)$.
     The image of $\pi^*$ is contained in the set of invariant functions:
     \[
     \image \pi^{*} \subseteq
     \C[X(\Gamma,\SL\C)]^{H^1(\Gamma,C_2)}.
     \]
     More precisely, we have $\pi^{*}f (\epsilon\chi)=
     f(\epsilon^{2}\chi^{2}) = \pi^{*}f (\chi)$ for all
     $\epsilon\in H^{1}(\Gamma,C_{2})$.
     It remains to prove that this inclusion is an equality.

     Since $ \C[X(\Gamma,\SL\C)]$ is
     generated  as $\C$-algebra by the functions $I_{\gamma}$ with
     $\gamma\in\Gamma$,  the monomials
     \[
         I_{\gamma_1} I_{\gamma_2}\cdots I_{\gamma_k}
     \]
     generate $\C[X(\Gamma,\SL\C)]$ as a $\C$-vector space.
     Taking the average of the action of $H^1(\Gamma,C_2)$, we deduce
that the subspace of invariant
     functions $\C[X(\Gamma,\SL\C)]^{H^1(\Gamma,C_2)}$ is generated by
     \[
     \frac 1 {2^r} \sum_{\epsilon\in H^1(\Gamma,C_2)}
            \epsilon \cdot I_{\gamma_{1}}\cdots
             I_{\gamma_{k}} =
             \big(\frac 1 {2^r}\sum_{\epsilon\in H^1(\Gamma,C_2)}
             \epsilon (\gamma_{1}\cdots\gamma_{k}) \big)
              I_{\gamma_{1}}\cdots
             I_{\gamma_{k}}
     \]
where $r$ is the rank of $H^{1}(\Gamma,C_2)$
         (see \cite[II.3.6]{Kraft} for instance).
         Using the fact that
         \[
         \frac 1 {2^r}\sum_{\epsilon\in H^1(\Gamma,C_2)}
         \epsilon(\gamma) =
         \begin{cases}
             1 & \text{ if } \gamma\in\Gamma^{2}\\
             0         & \text{ otherwise }
         \end{cases}
         \]
we deduce that $\C[X(\Gamma,\SL\C)]^{H^1(\Gamma,C_2)}$ is generated
by the monomials
     $
     I_{\gamma_1}\cdots I_{\gamma_k}
     $
such that the product $\gamma_1\ldots\gamma_k\in\Gamma^2$.

On the other hand we have for
$\chi\in X(\Gamma,\SL\C)$:
$$\pi^{*}\Sigma_{\gamma_{1},\ldots,\gamma_{k}} (\chi)
= \Sigma_{\gamma_{1},\ldots,\gamma_{k}}(\chi^{2}) =
\chi(\gamma_{1})\cdots\chi(\gamma_{k}) = I_{\gamma_{1}}\cdots
I_{\gamma_{k}} (\chi)\, ,$$ where
$\Sigma_{\gamma_{1},\ldots,\gamma_{k}}$ is the function defined in
(\ref{eqn:Sigma}). This gives that the monomials
$I_{\gamma_1}\cdots I_{\gamma_k}$ such that the product
$\gamma_1\ldots\gamma_k\in\Gamma^2$ is in the image 
of $\pi^{*}$
and therefore $\C[X(\Gamma,\SL\C)]^{H^1(\Gamma,C_2)}=\image
\pi^{*}$.
\end{proof}

\begin{rem}
Let $p\co X(\Gamma,\SL\C)\to\overline X(\Gamma)$ denote the
projection. If $\chi\in \overline X(\Gamma)$ is $\Ad$-irreducible,
then $p^{-1}(\chi)$ has $2^{r}$ points where $r$ is the rank of
$H^{1}(\Gamma,C_{2})$.
If $\chi$ is $\Ad$-reducible then the cardinality
of $p^{-1}(\chi)$ is strictly less than $2^{r}$. Thus $p$ is a
branched covering with branching locus the set of $\Ad$-reducible
characters.
\end{rem}

\begin{ex}
\label{ex:F2}    Let $F_2$ be the  free group of rank 2, with generators
$\alpha$ and $\beta$. There is an isomorphism:
\[
   (I_{\alpha},I_{\beta},I_{\alpha\beta})\co X(F_2,\SL\C)\to\C^3
\]
where $I_{\gamma}$ denotes the regular function induced by
$\tilde\tau_{\gamma}$.
In particular $X(F_2,\SL\C)$ is smooth.

Since every representation in $R(F_2)$ lifts to $\SL\C$, we deduce
\[
X(F_2)=X(F_2,\SL\C)/\!/{H^{1}(F_2,C_{2})}.
\]
The group $H^{1}(F_2,C_2)\cong ({C_2})^2$ has four elements, and its
action on $X(F_2,\SL\C)$ is generated by the involutions
   \[
\begin{array}{rcl}
(I_{\alpha},I_{\beta},I_{{\alpha}{\beta}})&\mapsto&(-I_{\alpha},I_{\beta},-I_{{\alpha}{\beta}})
\\
(I_{\alpha},I_{\beta},I_{{\alpha}{\beta}})&\mapsto&(I_{\alpha},-I_{\beta},-I_{{\alpha}{\beta}}).
\end{array}
   \] Thus $\C[X(F_2),\SL\C]^{H^{1}(F_2,C_{2})}$
is generated by $X=I_{\alpha}^2$, $Y=I_{\beta}^2$,
$Z=I_{{\alpha}{\beta}}^2$ and $W=I_{\alpha} I_{\beta}
I_{{\alpha}{\beta}}$. Hence
\begin{equation}
   \label{eqn:XF2}
   X(F_2)\cong\{(X,Y,Z,W)\in\C^4\mid W^2=XYZ\}
\end{equation}
The relationship with Corollary~\ref{cor:embedding}
   is given by the change of coordinates (cf.
Equality~(\ref{eqn:sigmaABAB}))
\[\left\{
\begin{array}{l}
      J_{\alpha}=X\\
      J_{\alpha}=Y\\
      J_{{\alpha}{\beta}}=Z\\
      J_{{\alpha}{\beta}^{-1}}=XY+Z-2W.
\end{array}
\right.
\]
\end{ex}

\begin{rem}
      From Equality~(\ref{eqn:XF2}) we remark that the singular
      set of $X(F_2 )$ consists of
      those points such that two of $\{ X,Y,Z\}$ vanish. This is the
      same as the set of characters of  representations generated by
      two rotations of angle $\pi$. This is also the set of
      $\Ad$-reducible but non-reducible representations.
\end{rem}

\begin{ex} If $M$ is a knot exterior in $S^3$, then $H_2(\pi_1M)\cong H_2(M)\cong
0$ and therefore $X(M)=\overline X(M)$. When in addition $M$ is a
 2-bridge knot exterior, explicit
methods of how to compute $X(M)$ are given in \cite{HLM1} and
\cite{HLM2}, where $X(M)$ for this particular case was already
defined as $X(M,\SL\C)/\!/C_2$. The explicit computation for the
figure eight knot exterior is found in \cite{GM}, for instance.
\end{ex}

\subsection{Representations that do not lift}

\begin{proof}[Proof of Theorem~\ref{thm:nolifts}]
The manifold
$M$ is a bundle over $S^1$ with fiber $\dot{T}^2$ a torus minus a
disk. Up to homeomorphism, $M$ is described by the action of the
monodromy on $H_1(\dot{T}^2,\mathbb Z)$, which is given by the
matrix
\[
\left(\begin{matrix} 1 & m_2 \\ m_1 & 1+m_1 m_2
\end{matrix}
\right)
\]
with $m_i\in 2\mathbb Z$, $m_i>0$. We shall show that
$X(M)-\overline X(M)$ has arbitrarily many components by choosing
$m_i$ sufficiently large.

   To have a presentation of $\pi_1M$, we use an automorphism $f$ of
$\pi_1\dot T^2$ induced by the monodromy. Since $\pi_1\dot T^2$ is
the free group of rank $2$ generated by $\alpha$ and $\beta$,
\[
\pi_1M=\langle \alpha,\beta,\mu\mid \mu\alpha\mu^{-1}=f(\alpha),
\mu\beta\mu^{-1}=f(\beta)\rangle
\]
We choose $f$ such that:
\[
\left\{
\begin{array}{rcl}
      \mu\alpha\mu^{-1}&=&\alpha\beta^{m_2}\\
      \mu\beta\mu^{-1}&=&\beta(\alpha\beta^{m_2})^{m_1}
\end{array}
\right.
\]
We choose odd numbers $p_1,p_2\in  2\mathbb Z+1$, with $1\leq
p_i\leq m_i/2$ and an arbitrary complex number $z\in \C$. By
Example~\ref{ex:F2}, there exist matrices $A_{z},B_{z}\in \SL\C$ with
\[
\tr(A_{z})=2\cos\frac{\pi p_1}{m_1},\
\tr(B_{z})=2\cos\frac{\pi p_2}{m_2} \text{ and }
\tr(A_{z} B_{z})=z\,.
\]
Those trace equalities imply that
$A_{z}^{m_1}=B_{z}^{m_2}=-\Id$. In particular
\[
\begin{array}{rcl}
      A_{z} B_{z}^{m_2} &=& -A_{z},\\
      B_{z} (A_{z} B_{z}^{m_2})^{m_1}&=&-B_{z}.
\end{array}
\]
Let $\rho_{z}\in R(\Gamma)$ be the representation that
$\rho_{z}(\alpha)=\pm A_{z} $, $\rho_z(\beta)=\pm B_{z}$ and
$\rho_z(\mu)=\pm \Id$. Since $m_1$ and $m_2$ are even, this
representation does not lift to $\SL\C$. In addition, for each
value of $p_1$ and $p_2$ we have defined a one parameter family of
characters, with parameter $z=\tr(A_{z}B_{z})\in\C$. By
\cite[Proposition~2.4]{CCGLS} the dimension of each component of
$X(M)$ is at most one, hence different values of $p_1$ and $p_2$
give different components.
\end{proof}

\section{The singular set of $X(F_n)$}

In this section we compute the singular set of $X(F_n)$, but
before we need two preliminary subsections: in
Subsection~\ref{ss:ZariskiTS} we recall some basic facts about the
Zariski tangent space and Luna's slice theorem, and in
Subsection~\ref{ss:cohomology} we compute the cohomology of free
groups with twisted coefficients.

\subsection{The Zariski tangent space}
\label{ss:ZariskiTS}

Given a representation $\rho\in R(\Gamma)$, we define the space of
cocycles
\[
      Z^1(\Gamma,\Ad\circ\rho) = \left\{
      \theta\co\Gamma\to \Sl\C\ \left\vert
      \begin{array}{c}
      \theta(\gamma_1\gamma_2)=\theta(\gamma_1)+\Ad_{\rho(\gamma_1)}
(\theta(\gamma_2)),
      \\
      \forall \gamma_1,\gamma_2\in\Gamma
      \end{array}
      \right.
      \right\}\, .
\]
Given a smooth path of representations $\rho_t$, with $t$ in a
neighborhood of the origin, one can construct a cocycle as
follows:
\[
\begin{array}{rcl}
          \Gamma &\to &\Sl\C\\
          \gamma&\mapsto&
          \frac{d\phantom{t}}{dt}\rho_t(\gamma)\rho_0(\gamma)^
{-1}\vert_{t=0}
\end{array} \,.
\]
This construction defines an isomorphism, due to Weil \cite{Weil}:

\begin{thm}[\cite{Weil}]\label{thm:Weil} The previous construction
defines an
isomorphism
\[
\TaZ\rho{R(\Gamma)}\cong Z^1(\Gamma,\Ad\circ\rho).
\]
\end{thm}

   Here $\TaZ\rho{R(\Gamma)}$ denotes the Zariski tangent space in the
scheme
sense (i.e.\ the defining ideals are not necessary reduced).

We also consider the space of coboundaries
\[
      B^1(\Gamma,\Ad\circ\rho) = \left\{
      \theta\co\Gamma\to\R^2\ \left\vert
      \begin{array}{c}
      \textrm{ there exists } a\in \Sl\C\textrm{ such that }\\
      \theta(\gamma)=\Ad_{\rho(\gamma)} (a)- a,\,
         \forall \gamma\in\Gamma
      \end{array}
      \right.
      \right\}\,.
\]

The isomorphism of Theorem~\ref{thm:Weil} identifies the subspace of the Zariski
tangent space corresponding to the orbits by conjugation with $
B^1(\Gamma,\Ad\circ\rho)$. 
So it seems natural that in some cases
$\TaZ\chi{X(\Gamma)}$ is isomorphic to the cohomology group
   \[
H^1(\Gamma,\Ad\circ\rho)=
Z^1(\Gamma,\Ad\circ\rho)/B^1(\Gamma,\Ad\circ\rho)
   \]
as we will show next.

The stabilizer of a representation $\rho\in R(\Gamma)$ is denoted
by
\[
Stab_{\rho}=\{A\in \PSL\C\mid A\rho A^{-1}=\rho\}\,.
\]
In particular, for and $\Ad$-irreducible representation
$Stab_{\rho}$ is trivial.

\begin{prop}
\label{prop:TZ}
   If $\rho$ is a smooth point of $R(\Gamma)$ with closed orbit,
then
\[
\TaZ{\chi_{\rho}}{X(\Gamma)}\cong
\TaZ0{H^1(\Gamma,\Ad\circ\rho)/\!/ Stab_{\rho}}\,.
\]
\end{prop}

\begin{proof} We use the slice theorem of Luna: there exists an
algebraic subvariety $S\subset R(\Gamma)$ that contains $\rho$ and
that is $  Stab_{\rho}$-invariant, such that
\begin{equation}
      \label{eqn:complement}
Z^1(\Gamma,\Ad\circ\rho)=B^1(\Gamma,\Ad\circ\rho)\oplus
\TaZ{\rho}{S}
\end{equation}
and the map induced by the projection
\[
S/\!/Stab_{\rho}\to X(\Gamma)
\]
is an \'etale isomorphism  
(in particular their tangent spaces are
isomorphic). Since we assume that $\rho$ is a smooth point, Luna's
theorem shows  that $S/\!/Stab_{\rho}$ and 
$\TaZ{\rho}{S}/\!/Stab_{\rho}$ are \'etale equivalent (see \cite[p. 97
]{KSS}). Since $\TaZ{\rho}{S}$ and $H^1(\Gamma,\Ad\circ\rho)$ are
isomorphic as $Stab_{\rho}$-modules (by
Equation~(\ref{eqn:complement})), the proposition follows.
\end{proof}

\subsection{Cohomology of Free groups}\label{ss:cohomology}

We start with irreducible characters:

\begin{lem}
\label{lem:cohomologyirreducible}
   Let $\chi_{\rho}\in X(F_n)$ be an irreducible
character. Then
\[
      \dim H^1(F_n,\Ad\circ\rho)=3n-3.
\]
\end{lem}

\begin{proof} Notice first that $Z^1(F_{n},\Ad\circ\rho)\cong
\Sl\C^n\cong \C^{3n}$. Irreducibility implies that $\dim
B^1(F_{n},\Ad\circ\rho)=3$, which is  maximal (even if $\Ad$-reducible
representations have invariant subspaces, irreducibility implies
that the eigenvalues are different from $1$). 
\end{proof}

We are interested in computing $H^1(F_n,\Ad\circ\rho)$ 
as a $Stab_{\rho}$-module. If $ \rho$ is  $\Ad$-irreducible, then
$Stab_{\rho}$ is trivial, and therefore $H^1(F_n,\Ad\circ\rho)$ is the
trivial module $\C^{3n-3}$. In the reducible and $\Ad$-reducible cases
we need further computations.

\medskip

\underline{Reducible characters}. Let $\chi\in X(F_n)$ be a non
trivial reducible character. There exists a representation
$\rho\in R(F_n)$  with character $\chi$ such that $\rho$ consists
of diagonal matrices, constructed in Remark~\ref{rem:rhodiagonal}.

   We decompose the Lie algebra $\Sl\C=h_0\oplus h_-\oplus h_+$,
   where $h_0$, $h_+$ and $h_-$ are the one dimensional $\C$-vector
   spaces generated respectively by $ \left(\!\begin{smallmatrix} 1&0\\
0&-1
   \end{smallmatrix}\!\right)$, $ \left(\!\begin{smallmatrix} 0&1\\
   0&0
   \end{smallmatrix}\!\right)$ and $ \left(\!\begin{smallmatrix} 0&0\\
   1&0
   \end{smallmatrix}\!\right)$.

\begin{lem}
\label{lem:splitdiagonal} If $\rho$ is diagonal then
$\Ad\circ\rho$ preserves the splitting $\Sl\C=h_0\oplus h_-\oplus
h_+$. If in addition $\rho$ is non-trivial, then $Stab_{\rho}$
preserves the splitting $\Sl\C=h_0\oplus (h_-\oplus h_+)$ (some
elements may permute $h_+$ and $h_-$).
\end{lem}

\begin{proof}
The first assertion is clear, because diagonal matrices preserve
each factor $h_0$ and $h_{\pm}$.

When the image of $\rho$ has order $\geq 3$, the group $Stab_{\rho}$
is precisely the set of diagonal matrices. When the image has order
precisely $2$, then $Stab_{\rho}$ is the group of diagonal and
anti-diagonal ones $\left(\!\begin{smallmatrix} 0&*\\ *&0
\end{smallmatrix}\!
\right)$. Antidiagonal matrices preserve $h_0$ and permute $h_-$
with $h_+$, hence the second assertion is proved.
\end{proof}

\begin{lem}
\label{lem:cohmologydiagonal}
   Let $\rho\in R(F_n)$ be a non-trivial diagonal representation, then
   $H^1(F_n,\Ad\circ\rho)\cong h_0^n\oplus (h_+\oplus h_-)^{n-1} $
as $Stab_{\rho}$-modules.
\end{lem}

\begin{proof} By construction,  $Z^1(F_n,\Ad\circ\rho) \cong \Sl\C^n$.
We have the splitting
   \[ H^1(F_n,\Ad\circ\rho)\cong H^1(F_n,h_0)\oplus
H^1(F_n,h_+)\oplus H^1(F_n,h_-).
\]
A diagonal matrix $\pm \left(\!\begin{smallmatrix} a&0 \\ 0&a^{-1}
\end{smallmatrix}\!\right)$ acts trivially on $h_0$ and by
multiplication by a factor $a^{\pm 2}$ on $h_{\pm}$. Therefore
$B^1(F_n,h_0)\cong 0$ and $B^1(F_n,h_{\pm})\cong h_{\pm}$, and the
lemma follows.
\end{proof}

\underline{$\Ad$-reducible but irreducible characters}. Let
$\rho\in R(\Gamma)$ be irreducible but $\Ad$-reducible. Up to
conjugation the image of $\rho$ is contained in the group of
diagonal and anti-diagonal matrices. There are two possibilities
for the stabilizer $Stab_{\rho}$. If the image of $\rho$ has more
than four elements, then $Stab_{\rho}$ has two elements: the
identity and $\pm\left(\!\begin{smallmatrix} i&0\\ 0&-i
\end{smallmatrix}\!\right)$. Otherwise the image of $\rho$ is
Klein's 4-group (i.e.\ the group generated by
$\pm\left(\!\begin{smallmatrix} i&0\\ 0&-i
\end{smallmatrix}\!\right)$ and $\pm\left(\!\begin{smallmatrix} 0&1\\
-1&0
\end{smallmatrix}\!\right)$
). In this case $Stab_{\rho}$ equals the image of $\rho$.
   With the same argument as in
Lemma~\ref{lem:splitdiagonal}, one can prove:

\begin{lem} Let $\rho$ be as above. Then both $\Ad\circ\rho$
and $Stab_{\rho}$ preserve  the splitting $\Sl\C=h_0\oplus
(h_+\oplus h_-)$.
\end{lem}

\begin{lem}
\label{lem:cohomologyantidiagonal} Let $\rho\in R(F_n)$ be an
irreducible but $\Ad$-reducible representation, then
$H^1(F_n,\Ad\circ\rho)\cong \Sl\C^{n-1}$ as $Stab_{\rho}$-modules.
\end{lem}

\begin{proof}
   Again  $Z^1(F_n,\Ad\circ\rho) \cong \Sl\C^n$,
and we have the decomposition \[ H^1(F_n,\Ad\circ\rho)\cong
H^1(F_n,h_0)\oplus H^1(F_n,h_+\oplus h_-).
\]
The group $B^1(F_n,h_0)$ has dimension one, because the
antidiagonal matrices act on $h_0$ by change of sign. In addition,
$\dim(B^1(F_n,h_+\oplus h_-))=2$ is also maximal, because this is
the case when we restrict it to diagonal representations (see the
proof of Lemma~\ref{lem:cohmologydiagonal}).
\end{proof}

\subsection{Singular locus for free groups}

We saw above that $X(F_2,\SL\C)\cong\C^3$  is smooth. We also showed
that the singular points of $X(F_2 )$ are $\Ad$-reducible but
irreducible characters.

\begin{prop}
For $n\geq 3$ the singular set of $X(F_n )$ is precisely the set of
$\Ad$-reducible characters.
\end{prop}

\begin{proof}
Since $R(F_n)\cong \PSL\C^{n}$, $X(F_n)$ is irreducible and of
dimension $3 n-3$. 
Thus  $\chi\in X(F_n)$ is singular if and only
if
\[
\dim T^{Zar}_{\chi} X(F_n)> 3 n -3.
\]
This dimension is computed by means of Proposition~\ref{prop:TZ}: if
the orbit of $\rho\in t^{-1}(\chi)$ is closed then
\[
\dim T^{Zar}_{\chi} X(F_n)=\dim
T_0^{Zar}(H^1(F_n,\Ad\circ\rho)/\!/Stab_{\rho}).
\]
If $\rho\in R(F_n)$ is irreducible, by
Lemma~\ref{lem:cohomologyirreducible}
   $\dim H^1(F_n,\Ad\circ\rho)=3n-3$. If in addition $\rho$ is
$\Ad$-irreducible, then $Stab_{\rho}$ is trivial and therefore
$\chi_{\rho}$ is smooth.

If $\rho$ is irreducible but  $\Ad$-reducible, then
$H^1(F_n,\Ad\circ\rho)\cong \Sl\C^{n-1} $ as $Stab_{\rho}$ modules,
by Lemma~\ref{lem:cohomologyantidiagonal}. We may assume that the
image of $\rho$ has more than 4 elements, because the adherence set of
such characters is the whole set of irreducible but
$\Ad$-reducible characters, and the singular set is closed. Hence
$Stab_{\rho}$ is the group generated by the involution $\pm
\left(\!\begin{smallmatrix} i&0\\ 0&-i
\end{smallmatrix}\!
\right)$, that acts trivially on $h_0$ but as a change of sign on
$h_+\oplus h_-$. Thus the action of $Stab_{\rho}$ on
$H^1(F_n,\Ad\circ\rho)$ is equivalent to the involution on $\C^{3n-3}$
that fixes $(n-1)$ coordinates and changes the sign of the remaining
$(2n-2)$ coordinates. The quotient of $\C^{3n-3}$ by this involution
is not smooth, hence \(\dim
T_0^{Zar}(H^1(F_n,\Ad\circ\rho)/\!/Stab_{\rho})>3n-3\).

When $\chi_{\rho}$ is reducible but non trivial, we may assume
that $\rho$ is diagonal and its image has more that three elements
(again the adherence set of those characters is the whole set of
reducible ones). Thus $Stab_{\rho}$ is   the group of diagonal
matrices, and by Lemma~\ref{lem:cohmologydiagonal},
$H^1(F_n,\Ad\circ\rho)\cong h_0^n\oplus (h_+\oplus h_-)^{n-1}$ as
$Stab_{\rho}$-module. We have an isomorphism $Stab_{\rho}\cong
\C^*$ and
   $t\in \C^*$
   acts on $h_{0}$ trivially and on $h_{\pm}$ by multiplication by
$t^{\pm 1}$.  An elementary computation shows that $(h_+\oplus
h_-)^{n-1}/\!/\C^*$ has dimension $2n-3$ and it is not smooth for
$n>2$.
\end{proof}

A similar argument yields that for $n\geq 3$ the singular part of
$X(F_n,\SL\C)$ is precisely the set of reducible characters.

\bigskip

\textbf{Acknowledgment} The second author was partially supported by
MCYT through grant BFM2000-0007.

\bigskip

\textsc{Laboratoire de Math\'ematiques Pures, Universit\'e Blaise
Pascal, F-63177, Aubi\`ere Cedex, France},
heusener@math.univ-bpclermont.fr.

\smallskip

\textsc{Departament de Matem\`atiques, Universitat Aut\`onoma de
Barcelona, 08193 Bellaterra, Spain},  porti@mat.uab.es.


\begin{thebibliography}{[CCGLS]}

\bibitem[ATZ]{ATZ}
H.\  Aslaksen, E.\  Tan, C.\  Zhu,
\newblock
Invariant theory of special orthogonal groups.
\newblock
\emph{Pacific J. Math.} \textbf{168} (1995),  207--215.


\bibitem[Bea]{Beardon}
A.F.\ Beardon.
\newblock \emph{The geometry of discrete groups}.
\newblock GTM 91, Springer Verlag New York, Berlin Heidelberg,
1983.


\bibitem[BMP]{BMP}
M.\ Boileau, S.\ Maillot, J.\ Porti.
\newblock
  \emph{Three-dimensional orbifolds and their geometric structures},
\newblock Soc.\ Math.\ France, Paris, 2003.
      


\bibitem[BZ]{BZ}
S.\ Boyer, X.\ Zhang,
\newblock
   On Culler-Shalen seminorms and Dehn filling.
\newblock
\emph{Ann. of Math. (2)} \textbf{148} (1998),
     737--801.

\bibitem[Bur90]{Bur90}
G.\ Burde.
\newblock
$SU(2)$-representation spaces for two-bridge knot groups.
\newblock \emph{Math. Ann.} \textbf{288} (1990),  103--119.

\bibitem[Cul]{Culler}
M.\ Culler.
\newblock Lifting representations to covering groups.
\newblock {\em Adv. in Math.} \textbf{59} (1986), 64--70.


\bibitem[CS]{CS}
M.\ Culler, P.B.\ Shalen.
\newblock Varieties of group representations and splittings of
3-manifolds.
\newblock {\em Ann. of Math. (2)} \textbf{117} (1983), 109--146.


\bibitem[CCGLS]{CCGLS}
D.~Cooper, M.~Culler, H.~Gillet, D.D.~Long,
P.B.~Shalen.
\newblock
   Plane curves associated to character varieties of
$3$-manifolds.
\newblock \emph{Invent. Math.} \textbf{118} (1994),  47--84.


\bibitem[FH]{FultonHarris}
W.\ Fulton, J.\ Harris.
\newblock \emph{Representation Theory}.
\newblock GTM 129, Springer Verlag New York, Berlin Heidelberg,
1991.



\bibitem[GM]{GM}
F.\ Gonz\'alez-Acu\~na, J.M.\ Montesinos-Amilibia.
\newblock On the character variety of group representations in
${\mathrm SL}(2,{\mathbf C})$ and ${\mathrm PSL}(2,{\mathbf C})$.
\newblock \emph{Math. Z.} \textbf{214} (1993), 627--652.


\bibitem[HLM1]{HLM1}
H.M.\ Hilden, M.T.\ Lozano, J.M.\ Montesinos-Amilibia.
\newblock On the character variety of group representations of a
$2$-bridge link $p/3$ into ${\rm PSL}(2, C)$. \newblock \emph{  Bol.
Soc. Mat. Mexicana}  \textbf{37}  (1992), 241--262.



\bibitem[HLM2]{HLM2}
H.M.\ Hilden, M.T.\ Lozano, J.M.\ Montesinos-Amilibia.
\newblock On the arithmetic $2$-bridge knots and link orbifolds and a new knot invariant.
  \newblock \emph{J. Knot Theory Ramifications}   \textbf{4}  (1995),
  81--114.



\bibitem[Kra]{Kraft}
H.\   Kraft.
\newblock
\emph{Geometrische Methoden in der Invariantentheorie.}
\newblock {Aspects of Mathematics.} Vieweg Verlag, Braunschweig,
     1985.


\bibitem[KSS]{KSS}
H.\   Kraft, P.\  Slodowy, T.A.\  Springer.
\newblock
\emph{Algebraische Transformationsgruppen und Invariantentheorie.}
\newblock {DMV Seminar,} \textbf{13}. Birkh\"auser Verlag, Basel,
     1989.

\bibitem[PV]{PV}
V.L.~Popov, \`E.B.~Vinberg. 
\newblock
\emph{Invariant theory}. 
\newblock In {\em Algebraic Geometry~IV.} 
\newblock{EMS, 55. Springer-Verlag, Berlin, 1994.}


\bibitem[Ril84]{Ril84}
R.\ Riley, \newblock{Non abelian representations of 2-bridge knot groups}.
\newblock{\emph{Quart. J. Math. Oxford} \textbf{2} (1984) 191--208.}



\bibitem[Sha]{Sha}
I.R. Shafarevich.
\newblock {\em Basic Algebraic Geometry}.
\newblock Springer Verlag, 1977.



\bibitem[Weil]{Weil}
A.\ Weil, \newblock{Remarks on the cohomology of groups}.
\newblock{\emph{Ann. of Math.} \textbf{80} (1964) 149--157.}

\end{thebibliography}
\end{document}